\documentclass[12pt,reqno]{amsart}
\usepackage{amssymb,verbatim,enumerate,ifthen,amsthm,fancyhdr,mathtools}
\usepackage[mathscr]{eucal}
\usepackage[utf8]{inputenc}
\usepackage[T1]{fontenc}
\def\N{\mathbb{N}}
\def\R{\mathbb{R}}
\def\Q{\mathbb{Q}}
\def\Z{\mathbb{Z}}

\def\A{\mathscr{A}}

\newtheorem{theorem}{Theorem}
\newtheorem*{theorem*}{Theorem}
\def\Thm#1#2{\ifthenelse{\equal{#1}{*}}{\begin{theorem*}#2\end{theorem*}}
             {\begin{theorem}\label{T#1}#2\end{theorem}}}
\newtheorem{Atheorem}{Theorem}

\def\thm#1{Theorem~\ref{T#1}}

\newtheorem{proposition}[theorem]{Proposition}
\newtheorem*{proposition*}{Proposition}
\def\Prp#1#2{\ifthenelse{\equal{#1}{*}}{\begin{proposition*}#2\end{proposition*}}
             {\begin{proposition}\label{P#1}#2\end{proposition}}}

\newtheorem{corollary}[theorem]{Corollary}
\newtheorem*{corollary*}{Corollary}
\def\Cor#1#2{\ifthenelse{\equal{#1}{*}}{\begin{corollary*}#2\end{corollary*}}
             {\begin{corollary}\label{C#1}#2\end{corollary}}}
\def\cor#1{Corollary~\ref{C#1}}

\newtheorem{lemma}[theorem]{Lemma}
\newtheorem*{lemma*}{Lemma}
\def\Lem#1#2{\ifthenelse{\equal{#1}{*}}{\begin{lemma*}#2\end{lemma*}}
             {\begin{lemma}\label{L#1}#2\end{lemma}}}
\def\lem#1{Lemma~\ref{L#1}}
\newtheorem{Alemma}{Lemma}

\theoremstyle{definition}
\newtheorem{remark}[theorem]{Remark}
\newtheorem*{remark*}{Remark}
\def\Rem#1#2{\ifthenelse{\equal{#1}{*}}{\begin{remark*}\rm #2\end{remark*}}
             {\begin{remark}\label{R#1}\rm #2\end{remark}}}

\newtheorem{example}[theorem]{Example}
\newtheorem*{example*}{Example}
\def\Exa#1#2{\ifthenelse{\equal{#1}{*}}{\begin{example*}\rm #2\end{example*}}
             {\begin{example}\label{Ex#1}\rm #2\end{example}}}

\def\eq#1{{\rm(\ref{E#1})}}
\def\Eq#1#2{\ifthenelse{\equal{#1}{*}}
  {\begin{equation*}\begin{aligned}#2\end{aligned}\end{equation*}}
  {\begin{equation}\begin{aligned}\label{E#1}#2\end{aligned}\end{equation}}}

\def\dom{\mathop{\hbox{\rm dom}}\nolimits}

\begin{document}
\vspace{5mm}

\date{\today}

\title[On additive functions]{On additive functions with additional derivation properties}

\author[R.\ Gr\"unwald]{Rich\'ard Gr\"unwald}
\author[Zs. P\'ales]{Zsolt P\'ales}
\address{Institute of Mathematics, University of Debrecen, 
H-4002 Debrecen, Pf. 400, Hungary}
\email{richard.grunwald96@gmail.com, pales@science.unideb.hu}

\thanks{The research of the first author was supported by the \'UNKP-19-2 New National Excellence Program of the Ministry of Human Capacities. The research of the second author was supported by the 2019-2.1.11-T\'ET-2019-00049, the EFOP-3.6.1-16-2016-00022 and the EFOP-3.6.2-16-2017-00015 projects. The last two projects are co-financed by the European Union and the European Social Fund.}
\subjclass[2010]{Primary 39B22, 39B40, 39B50}
\keywords{algebraic derivation; derivation for trigonometric functions; derivation for hyperbolic functions}

\begin{abstract}
The purpose of this paper is to introduce the notion of a generalized derivation which derivates a prescribed family of smooth vector-valued functions of several variables. The basic calculus rules are established and then a result derived which shows that if a function $f$ satisfies an addition theorem whose determining operation is derivable with respect to an additive function $d$, then the function $f$ is itself derivable with respect to $d$. As an application of this approach, new proof of a generalization of a recent result of Maksa is obtained. We also extend the result of Nishiyama and Horinouchi and formulate two open problems.
\end{abstract}

\maketitle

\section{Introduction}

Derivations are additive and Leibniz-type mappings of a ring into itself. More 
precisely, if $(R,+,\cdot)$ is a ring, then a function $d:R\to R$ is called a \emph{derivation} if, 
for all $x,y\in R$,
\begin{eqnarray}
d(x+y)&=&d(x)+d(y), \label{ES} \\
d(x\cdot y)&=&d(x)\cdot y+x\cdot d(y). \label{EP}
\end{eqnarray}

Derivations are used in many branches of analysis and algebra. For instance, nonnegative information functions are constructed via real derivations (see  
Dar\'oczy--Maksa \cite{DarMak79}, Maksa \cite{Mak76}). Nonconstant functions that are convex with respect to families of power means are also obtained in terms of real derivations (see Maksa--P\'ales \cite{MakPal15}). Derivations are used to express the general solutions of certain functional equations (see Fechner--Gselmann \cite{FecGse12}, Gselmann \cite{Gse12}, \cite{Gse13a}, Halter-Koch \cite{Hal00}, \cite{Hal00b}, Jurkat \cite{Jur65}). Generalizations, such as higher-order derivations, bi-derivations and approximate or near-derivations were studied by Badora \cite{Bad06c}, Gselmann \cite{Gse14b}, Gselmann--P\'ales \cite{GsePal16}, and Maksa \cite{Mak81b}, \cite{Mak87c}. 

We say that a function $d:\R\to\R$ derivates a differentiable function $f:I\to\R$ if the functional equation 
\Eq{*}{
  d(f(x))=f'(x)d(x) \qquad(x\in I)
}
holds. In the pioneering papers \cite{Kur64,Kur65} Kurepa proved that if $d$ is an additive functions which derivates one of the maps $x\mapsto x^2$ or $x\mapsto x^{-1}$, then it satisfies the Leibniz Rule, i.e., it is a standard derivation.
This result was then extended by Nishiyama and Horinouchi in \cite{NisHor68},
who proved an analogous statement about the derivability of the power function $x\mapsto x^r$ with rational exponent different from $0$ and $1$. Boros and Erdei in \cite{BorErd05} proved that those additive functions that derivate the map $x\mapsto\sqrt{1-x^2}$, that is, satisfy the identity
\Eq{bor}{
		d(\sqrt{1-x^2})=-\frac{x}{\sqrt{1-x^2}}d(x) \qquad(x\in]-1,1[\,),
}
are also standard derivations. Maksa in \cite{Mak13a} showed that if an additive function derivates any of the exponential, hyperbolic or trigonometric functions, then is has to be a standard derivation again. A counterpart of this result was obtained by Grünwald and Páles in \cite{GruPal18}, where an analogous statement 
was established assuming Leibniz property instead of additivity.

The purpose of this paper is to introduce the notion of generalized derivation which derivates a prescribed family of smooth vector-valued functions of several variables. After establishing the basic calculus rules in \thm{inv}, we derive in \cor{add} a result which shows that if a function $f$ satisfies an addition theorem whose determining operation is derivable with respect to an additive function $d$, then the function $f$ is itself derivable with respect to $d$.  
Using this result, we will be able to give a completely new proof for the aforementioned result of Maksa \cite{Mak13a}. In addition, we also generalize this result, because we require the derivability of the exponential, hyperbolic or trigonometric functions only on small intervals. In the last section of our paper, we also offer a generalization of Nishiyama and Horinouchi by replacing power functions of the form $P\circ Q^{-1}$, where $P$ and $Q$ are polynomials with rational coefficients. Finally, we formulate some open problems.

\section{Generalized derivations and their properties}

For fixed $n,m\in\N$, the class of $n$-variable $\R^m$-valued \emph{admissible functions} is defined as follows:
	\Eq{*}{
		\A_n^m:=\{f:\Omega\to\R^m \mid \emptyset\neq\Omega\subset\R^n \mbox{ is open and }
		f\mbox{ is Fréchet differentiable on }\Omega\}
	}
	and we set
	\Eq{*}{
		\A  := \bigcup_{n,m=1}^\infty \A_n^m.
	}
The set $\Omega$ related to $f$ will be called the \emph{domain of $f$} and denoted by $\dom_f$. In general, for a vector $x\in\R^n$, we will denote the $i$th coordinate of $x$ by $x_i$, and for a function $f\in\A_n^m$, $f_j$ will stand for the $j$th coordinate function of $f$.

We say that a function $d:\R\to\R$ is a \emph{derivation with 
respect to an admissible function $f\in\A_n^m$} if, for all $x\in\dom_f$ and $j\in\{1,\dots,m\}$,
	\Eq{df+}{
		d\big(f_j(x)\big)=\partial_1 f_j(x)d(x_1)+\dots+\partial_n f_j(x)d(x_n)
	}
	holds. Furthermore we say that \emph{$d$ is a derivation with respect to $A\subseteq\A$}, if $d$ is a derivation with respect to each member of $A$. For any $n\in\N$ and $x\in\R^n$, define $d(x)$ by 
	\Eq{*}{
	d(x):=(d(x_1),\dots,d(x_n)).
	}
	Then \eq{df+} can simply be rewritten as
	\Eq{*}{
	  d(f(x))=f'(x)d(x),
	}
	where $f\in\A_n^m$ and $f'(x)$ denotes the Fréchet derivative of $f$ at $x$, which is an $n\times m$ matrix whose entries are the partial derivatives $\partial_i f_j$ at $x$.

One can immediately see that a function $d:\R\to\R$ is a standard derivation if and only if it is a derivation with respect to $S_2$ and $P_2$, where 
	\Eq{*}{
		S_2(x_1,x_2):=x_1+x_2 \qquad\mbox{and}\qquad P_2(x_1,x_2):=x_1x_2 \qquad((x_1,x_2)\in\R^2).
	}
In what follows, we will prefer the terminology $d$ is additive (resp.\ of Leibniz-type) whenever $d$ is a derivation with respect to $S_2$ (resp.\ $P_2$).


The following result, which is a significant extension of \cite[Lemma A]{GruPal18} collects the basic rules for functions that are derivable with respect to a fixed real function. In particular, its second assertion will be very useful for our purposes.

\Thm{inv}{
For any function $d:\R\to\R$, we have the following three assertions.
\begin{enumerate}[(i)]
	\item Let $n,m,k\in\N$ and $f\in\A_n^m$, $g\in\R_m^k$. If $d$ is a derivation with respect to $f$ and $g$, then $d$ is also a derivation with respect to $g\circ f$.
	\item Let $n,m,k\in\N$ and $f\in\A_n^m$, $g\in\R_m^k$ such that $f(\dom_f)$ is open. If $d$ is a derivation with respect to $f$ and $g\circ f$, then $d$ is also a derivation with respect to $g$ on $f(\dom_f)\cap\dom_g$.
	\item Let $n\in\N$ and $f\in\A_n^n$ with a continuous nowhere singular derivative. If $d$ is a derivation with respect to $f$, then $d$ is also a derivation with respect to its inverse $f^{-1}$.
\end{enumerate}
}

\begin{proof}
By the assumptions of (i), for all $x\in\dom_f$ and $y\in\dom_g$, we have
\Eq{xy}{
	d(f(x))=f'(x)d(x)\qquad\mbox{and}\qquad d(g(y))=g'(y)d(y).
}
Let $x\in\dom_f$ with $y:=f(x)\in\dom_g$. Using that $d$ is a derivation with respect to $f$ and $g$, by the standard Chain Rule, we get 
\Eq{*}{
	d((g\circ f)(x))&=d(g(f(x)))=d(g(y))=g'(y)d(y)\\
	&=g'(f(x))d(f(x))=g'(f(x))f'(x)d(x)=(g\circ f)'(x)d(x),
}
which yields that $d$ is a derivation with respect to the function $g\circ f$.

Let $y\in f(\dom_f)\cap\dom_g$. Then there exists $x\in\dom_f$ such that $y=f(x)$. Thus, applying the standard Chain Rule, we get
\Eq{*}{
  d(g(y))&=d(g(f(x))=d(g\circ f(x))=(g\circ f)'(x)d(x)\\
  &=g'(f(x))f'(x)d(x)=g'(f(x))d(f(x))=g'(y)d(y),
}
which proves that $d$ is also a derivation with respect to $g$ on $f(\dom_f)\cap\dom_g$.

By the assumption of (iii), for all $x\in\dom_f$, we have the first equality in \eq{xy}. Let $y\in\dom_{f^{-1}}$. 
Using the substitution $x=f^{-1}(y)$, this implies
\Eq{*}{
  d(y)=f'(f^{-1}(y))d(f^{-1}(y)).
}
Thus, by the inverse function theorem, it follows that
\Eq{*}{
	d(f^{-1}(y))=\big(f'(f^{-1}(y))\big)^{-1}d(y)=\big(f^{-1}\big)'(y)d(y).
} 
Thus, $d$ is a derivation with respect to the inverse function $f^{-1}$.
\end{proof}

The following consequence of the above theorem will be useful in several proofs.

\Cor{add}{Let $\Omega_1,\Omega_2\subseteq\R^n$ be nonempty open sets, let $f:\Omega_1\cup\Omega_2\cup(\Omega_1+\Omega_2)\to\R^m$ be a Fréchet differentiable function such that $f(\Omega_1)$ and $f(\Omega_2)$ are open. Assume that there exists a Fréchet differentiable function $g:f(\Omega_1)\times f(\Omega_2)\to\R^m$ such that $f$ satisfies the functional equation
	\Eq{addition}{
		f(x+y)=g(f(x),f(y)) \qquad ((x,y)\in\Omega_1\times\Omega_2).
	}
	Let $d:\R\to\R$ be an additive function which is a derivation with respect to $f$. Then $d$ is also a derivation with respect to $g$ on $f(\Omega_1)\times f(\Omega_2)$.}

\begin{proof} Assume that $d$ is an additive derivation with respect to $f$. By the additivity of $d$, we have that $d$ is a derivation with respect to the mapping
\Eq{*}{
  \Omega_1\times\Omega_2\ni(x,y)\mapsto f(x+y). 
}
Thus, the equality \eq{addition} implies that $d$ is a derivation with respect to the composition 
\Eq{*}{
  \Omega_1\times\Omega_2\ni(x,y)\mapsto g(f(x),f(y)).
}
On the other hand, $d$ is trivially a derivation 
with respect to the mapping
\Eq{*}{
  \Omega_1\times\Omega_2\ni(x,y)\mapsto(f(x),f(y)).
}
Applying the second assertion of the previous theorem, now it follows that $d$ is also a derivation with respect to $g$ on $f(\Omega_1)\times f(\Omega_2)$.
\end{proof}

\section{Localization theorems}

In the sequel, for a number $r\in\Q$, let $D_r$ denote the domain of the power function $x\mapsto x^r$, which is defined in the following way: If $r=m/n$, where  
$n\in\N$, $m\in\Z$ and $n,m$ are coprime, then let
\Eq{*}{
  D_r:=\begin{cases}
       \R \qquad&\mbox{if $n$ is odd and $m\geq0$},\\
       \R\setminus\{0\} \qquad&\mbox{if $n$ is odd and $m<0$},\\
       [0,\infty[ \qquad&\mbox{if $n$ is even and $m\geq0$},\\
       ]0,\infty[ \qquad&\mbox{if $n$ is even and $m<0$}.\\
       \end{cases}
}
A function $d:\R\to\R$ is said to be \emph{$\Q$-homogeneous} if, for all $x\in\R$ and $r\in\Q$, the equality $d(rx)=rd(x)$ holds. It is well-known that every additive function is automatically $\Q$-homogeneous.

\Lem{Nisalt}{
	Let $r\in\Q$, let $I\subseteq D_{r-1}$ be a nonempty open subset and $d:\R\to\R$ be a $\Q$-homogeneous function. Suppose that the equality
	\Eq{xr}{
		d(x^r)=rx^{r-1}d(x)
	}
	holds for all $x\in I$. Then it is also valid for all $x\in D_{r-1}$.
} 

\begin{proof}
	Assume that \eq{xr} holds for all $x\in I$. 
	Let $x\in D_{r-1}$ be arbitrary. If $x=0$, then $r\geq1$ and hence $x^r=0^r=0$. Thus, by $d(0)=0$, \eq{xr} is trivially valid. If $r=1$ So may assume that $x$ is an arbitrarily fixed element from $D_r\setminus\{0\}$. 
	Suppose that $r$ is of the form $r=m/n$ for some $n\in\N$ and $m\in\Z$. By the density of $\Q$ in $\R$, it is clear that the set
	\Eq{*}{
		\{q\in\Q\mid q\neq0,\, q^nx\in I\}
	}
	is nonempty. Let $q$ be a fixed element from it. Thus, using the validity of equation \eq{xr} on the interval $I$ and the $\Q$-homogeneity of $d$, we obtain
	\Eq{*}{
		q^md(x^r)=d(q^mx^r)=d((q^nx)^r)=r(q^nx)^{r-1}d(q^nx)=q^mrx^{r-1}d(x),
	}
	which simplifies to $d(x^r)=rx^{r-1}d(x)$. This is exactly the desired equality since $x$ was an arbitrary element from $D_r\setminus\{0\}$.
\end{proof}

The following result essentially was proved by Nishiyama and Horinouchi \cite{NisHor68}. The result concerning the particular cases $r=-1$ and $r=2$ were discovered by Kurepa in \cite{Kur64} and \cite{Kur65}.

\Lem{Nis}{
	Let $d:\R\to\R$ be an additive function, let $r\in\Q\setminus\{0,1\}$ and let $I\subseteq D_r$ be a nonempty open subinterval. Then $d$ is a standard derivation if and only if \eq{xr} is valid for all $x\in I$.
}

\begin{proof} Assume first that $d$ is a standard derivation. Then by an easy argument, it follows that \eq{xr} is valid for all $x\in D_r$.

Conversely, if \eq{xr} is valid for all $x\in I$, then, by \lem{Nisalt}, we get that \eq{xr} is valid for all $x\in D_r$. Now, the result of Nishiyama and Horinouchi \cite{NisHor68} implies that $d$ must be a standard derivation.
\end{proof}

\Lem{LRext}{
	Let $U\subseteq\R^2$ be a nonempty open subset and let $d:\R\to\R$ be a $\Q$-homogeneous function which satisfies the functional equation \eqref{EP}
	for all $(x,y)\in U$. Then \eqref{EP} also holds for all $x,y\in\R$.  
}

\begin{proof}
	Let $x,y\in\R$. Using the density of $\Q$ in $\R$, it is clear that there exist $p,q\in\Q\setminus\{0\}$ such that $(px,qy)\in U$. Applying equation \eq{P} for $px$ and $qy$, taking into consideration the $\Q$-homogeneity of $d$, we obtain that 
	\Eq{*}{
		pqd(xy)=d((px)(qy))=d(px)qy+pxd(qy)=pd(x)qy+pxqd(y).
	}
	Dividing by $pq$ the above equality, we get the statement.
\end{proof}

To prove our main result, which will extend the theorem of Maksa \cite{Mak13a}, the following lemma will also be needed.

\Lem{Ratio}{
	The sets
	\Eq{*}{
		U&:=\{x\in\R: x,\sqrt{1+x^2}\in\Q\},\\
		V&:=\{x\in\,]-\infty,-1[\,\cup\,]1,\infty[\,: x,\sqrt{x^2-1}\in\Q\},\\
		W&:=\{x\in\,]-1,1[\,: x,\sqrt{1-x^2}\in\Q\}
	}
	are dense in $\R$, in $\,]-\infty,-1[\,\cup\,]1,\infty[\,$, and in $\,]-1,1[\,$, respectively.
}

\begin{proof}
	To prove the density of $U$, let $x\in \,]-1,1[\,$, let $0<\varepsilon<\min(1+x,1-x)$ be arbitrary and denote
	\Eq{*}{
		I:=\left]\sqrt{\frac{1-x-\varepsilon}{1+x+\varepsilon}},\sqrt{\frac{1-x+\varepsilon}{1+x-\varepsilon}}\right[.
	}
	We are going to show that 
	\Eq{rr}{
		\frac{1-r^2}{1+r^2}\in \,]x-\varepsilon,x+\varepsilon[\,\cap\, U
		\qquad \mbox{for all}\qquad r\in \Q\cap I.
	}
	Indeed, let $r\in\Q\cap I$ be arbitrary. Then there exists $(m,n)\in\Z\times\N$ such that $r=\frac{m}{n}$. Using that $r$ is bounded by the endpoints of $I$, we easily get that
	\Eq{*}{
		s:=\frac{1-r^2}{1+r^2}\in\,]x-\varepsilon,x+\varepsilon[\,.
	}
	On the other hand, $s\in\Q$ and
	\Eq{*}{
		\sqrt{1-s^2}
		=\sqrt{1-\bigg(\frac{1-r^2}{1+r^2}\bigg)^2}
		=\sqrt{1-\bigg(\frac{n^2-m^2}{n^2+m^2}\bigg)^2}
		=\frac{2nm}{n^2+m^2}\in\Q,
	}
	which completes the proof of \eq{rr}.
	
	To prove the density of $V$, let $x\in \,]-\infty,-1[\,\cup\,]1,\infty[\,$, that is, let $|x|>1$, let $0<\varepsilon<|x|-1$ be arbitrary and denote
	\Eq{*}{
		J:=
		\begin{cases}
			\left]\sqrt{\dfrac{x+\varepsilon-1}{x+\varepsilon+1}},\sqrt{\dfrac{x-\varepsilon-1}{x-\varepsilon+1}}\right[ 
			&\mbox{if } x<-1,\\[4mm]
			\left]\sqrt{\dfrac{x-\varepsilon-1}{x-\varepsilon+1}},\sqrt{\dfrac{x+\varepsilon-1}{x+\varepsilon+1}}\right[ 
			&\mbox{if } 1<x.\\
		\end{cases}
	}
	Then, one can easily check that $J$ is nonempty and $J\subseteq\,]1,\infty[\,$ if $x<-1$ and $J\subseteq\,]0,1[\,$ if $1<x$. Thus $1\not\in J$ holds in both cases.
	We are going to show that 
	\Eq{rr2}{
		\frac{1+r^2}{1-r^2}\in \,]x-\varepsilon,x+\varepsilon[\,\cap\, V
		\qquad \mbox{for all}\qquad r\in \Q\cap J.
	}
	Indeed, let $r\in\Q\,\cap J$ be arbitrary. Then there exists $(m,n)\in\N\times\N$ such that $r=\frac{m}{n}$. Since $r$ cannot be equal to $1$, therefore $n\neq m$ holds. Using that $r$ is bounded by the endpoints of $J$, in each cases we easily get that
	\Eq{*}{
		s:=\frac{1+r^2}{1-r^2}\in\,]x-\varepsilon,x+\varepsilon[\,.
	}
	On the other hand, $s\in\Q$ and
	\Eq{*}{
		\sqrt{s^2-1}
		=\sqrt{\bigg(\frac{1+r^2}{1-r^2}\bigg)^2-1}
		=\sqrt{\bigg(\frac{n^2+m^2}{n^2-m^2}\bigg)^2-1}
		=\frac{2nm}{|n^2-m^2|}\in\Q,
	}
	which completes the proof of \eq{rr2}. 
	
	To prove the density of $W$, let $x\in\R\setminus\{0\}$ and $0<\varepsilon<|x|$ be arbitrary and denote
	\Eq{*}{
		K:=
		\begin{cases}
			\left]\dfrac{-1}{x+\varepsilon}-\sqrt{\dfrac{1}{(x+\varepsilon)^2}+1},1\right[ 
			&\mbox{if } x<0,\\[4mm]
			\left]-1,\dfrac{-1}{x+\varepsilon}+\sqrt{\dfrac{1}{(x+\varepsilon)^2}+1}\right[ 
			&\mbox{if } 0<x.\\
		\end{cases}
	}
	Then, one can easily check that $K$ is nonempty and $K\subseteq\,]-1,1[\,$. Thus $1\not\in K$ holds in both cases.
	We are going to show that 
	\Eq{rr3}{
		\frac{2r}{1-r^2}\in \,]x-\varepsilon,x+\varepsilon[\,\cap\, W
		\qquad \mbox{for all}\qquad r\in \Q\cap K.
	}
	Indeed, let $r\in\Q\cap K$ be arbitrary. Then there exists $(m,n)\in\Z\times\N$ such that $r=\frac{m}{n}$. Since $r$ cannot be equal to $1$, therefore $n\neq m$ holds. Using that $r$ is bounded by the endpoints of $K$, in each cases we easily get that
	\Eq{*}{
		s:=\frac{2r}{1-r^2}\in\,]x-\varepsilon,x+\varepsilon[\,.
	}
	On the other hand, $s\in\Q$ and
	\Eq{*}{
		\sqrt{1+s^2}
		=\sqrt{1+\bigg(\frac{2r}{1-r^2}\bigg)^2}
		=\sqrt{1+\bigg(\frac{2nm}{n^2-m^2}\bigg)^2}
		=\frac{n^2+m^2}{|n^2-m^2|}\in\Q,
	}
	which completes the proof of \eq{rr3}.
\end{proof}

\section{Extension of the result of Maksa}

In what follows we extend the result of Maksa \cite{Mak13a} by assuming the derivability of any of the exponential, hyperbolic or trigonometric functions on a small interval. Our approach is based on the use of the addition theorems for each of these functions and the application \cor{add}. In each particular case, we obtain that $d$ is a derivation with respect to a two-variable algebraic function.

\Thm{Mak}{
	Let $d:\R\to\R$ be an additive function and let $\alpha,\beta\in\R$ with $\alpha<\beta$. Suppose that $d$ is a derivation with respect to the restriction to $\,]\alpha,\beta[\,$ of any of the following functions with further assumptions on $\alpha$ and $\beta$, respectively:
	\Eq{Mak}{
	(i)\qquad\qquad\qquad&\exp \qquad&&(2\alpha<\beta\mbox{ and }\alpha<2\beta),\\
	(ii)\qquad\qquad\qquad	&\sinh \qquad&&(\alpha<0<\beta),\\
	(iii)\qquad\qquad\qquad	&\cosh \qquad&&(0<2\alpha<\beta\mbox{ or }\alpha<2\beta<0),\\
	(iv)\qquad\qquad\qquad	&\tanh \qquad&&(\alpha<0<\beta),\\
	(v)\qquad\qquad\qquad	&\coth \qquad&&(0<2\alpha<\beta\mbox{ or }\alpha<2\beta<0),\\
	(vi)\qquad\qquad\qquad	&\sin \qquad&&(\alpha<0<\beta),\\
	(vii)\qquad\qquad\qquad	&\cos \qquad&&(0<2\alpha<\pi<\beta\mbox{ or } \alpha<-\pi<2\beta<0),\\
	(viii)\qquad\qquad\qquad	&\tan \qquad&&(-\pi<2\alpha<\beta\mbox{ and }\alpha<2\beta<\pi),\\
	(ix)\qquad\qquad\qquad	&\cot \qquad&&(0<2\alpha<\beta<\pi\mbox{ or }-\pi<\alpha<2\beta<0).
	}
	Then $d$ is a standard derivation. 
}

\begin{proof} Observe that in each of the above cases the inequalities $2\alpha<\beta$ and $\alpha<2\beta$ hold. Adding up these inequalities side by side, it follows that $\alpha<\beta$ and hence
	\Eq{*}{
		\gamma:=\tfrac12\max(\alpha,2\alpha)
		<\tfrac12\min(\beta,2\beta)=:\delta.
	}
	Then $\,]\gamma,\delta[\,\subseteq \,]\alpha,\beta[\,\cap \,]\tfrac\alpha2,\tfrac\beta2[\,$, which implies that 
	\Eq{*}{
		\,]\gamma,\delta[\,+\,]\gamma,\delta[\,
		=\,]2\gamma,2\delta[\,\subseteq\,]\alpha,\beta[\,.
	}
	In the rest of proof, we shall utilize that each of the functions listed in \eq{Mak} possesses an addition formula, i.e., it satisfies  functional equation of type \eq{addition} with $\Omega_1:=\Omega_2:=\,]\gamma,\delta[\,$.
	
	$(i)$ Assume first that $2\alpha<\beta$ and $\alpha<2\beta$ and $d$ is a derivation with respect to the restriction to $]\alpha,\beta[\,$ of the exponential function. It means that the equation $d(\exp(x))=\exp(x)d(x)$ holds for all $x\in \,]\alpha,\beta[\,$. Using that this restriction satisfies the functional equation 
	\Eq{*}{
		\exp(x+y)=\exp(x)\exp(y) \qquad(x,y\in\,]\gamma,\delta[\,),
	} 
	\cor{add} implies that $d$ is a derivation with respect to the mapping 
	\Eq{*}{
		(u,v)\mapsto u\cdot v\qquad(u,v\in\,]\exp(\gamma),\exp(\delta)[\,),
	}
	i.e., $d$ is of Leibniz-type on the interval $]\exp(\gamma),\exp(\delta)[\,$. In view of \lem{LRext}, it follows that $d$ is of  Leibniz-type on $\R$ and hence it is a standard derivation.
	
	$(ii)$ In the second case assume that $\alpha<0<\beta$ and $d(\sinh(x))=\cosh(x)d(x)$ holds for all $x\in\,]\alpha,\beta[\,$. Then $\gamma=\frac\alpha2$ and $\delta=\frac\beta2$. First we can choose $\lambda>0$ such that $\,]-\lambda,\lambda[\,\subseteq \,]\gamma,\delta[\,$. Using the identity $\cosh(x)=\sqrt{1+\sinh^2(x)}$ ($x\in\,]-\lambda,\lambda[\,$), we obtain that the restriction of the sine hyperbolic function to the interval $\,]\alpha,\beta[\,$ satisfies the functional equation 
	\Eq{*}{
		\sinh(x+y)=\sinh(x)\sqrt{1+\sinh^2(y)}+\sqrt{1+\sinh^2(x)}\sinh(y)
	}
	for all $x,y\in\,]-\lambda,\lambda[\,$. By \cor{add}, it follows that $d$ is also a derivation with respect to the function 
	\Eq{*}{
		(u,v)\mapsto u\sqrt{1+v^2}+v\sqrt{1+u^2}\qquad (u,v\in\,]-\sinh(\lambda),\sinh(\lambda)[\,).
	}
	It means that the functional equation
	\Eq{Phisinh}{
		d(u&\sqrt{1+v^2}+v\sqrt{1+u^2}) \\&
		=\Big(\sqrt{1+v^2}+v\frac{u}{\sqrt{1+u^2}}\Big)d(u)+\Big(\sqrt{1+u^2}+u\frac{v}{\sqrt{1+v^2}}\Big)d(v)
	}
	holds for all $u,v\in\,]-\sinh(\lambda),\sinh(\lambda)[\,$. Replacing $v$ by $-v$ and adding the equality so obtained to \eq{Phisinh} side by side, we get that
	\Eq{Phisinh+}{
		d(u\sqrt{1+v^2}) 
		=\sqrt{1+v^2}d(u)+u\frac{v}{\sqrt{1+v^2}}d(v)
	}
	holds for all $u,v\in\,]-\sinh(\lambda),\sinh(\lambda)[\,$. Let $U$ be the set defined in \lem{Ratio}. Then, by this lemma, the intersection $U\cap \,]-\sinh(\lambda),\sinh(\lambda)[\,$ is nonempty, moreover it is dense in $]-\sinh(\lambda),\sinh(\lambda)[\,$. Then, for $u,v\in U\cap \,]-\sinh(\lambda),\sinh(\lambda)[\,$, we have that $u,v,$ and $\sqrt{1+v^2}\in\Q$. For such values of $u$ and $v$, the equality \eq{Phisinh+} and the $\Q$-homogeneity of $d$ implies that
	\Eq{*}{
		u\sqrt{1+v^2}d(1)
		=\sqrt{1+v^2}ud(1)+u\frac{v}{\sqrt{1+v^2}}vd(1),
	}
	which is possible only if $d(1)=0$. Let $0\neq u\in\,]-\sinh(\lambda),\sinh(\lambda)[\,\cap\,\Q$,
	and $v\in\,]-\sinh(\lambda),\sinh(\lambda)[\,$. By the $\Q$-homogeneity and additivity of $d$, it follows that $d(u)=ud(1)=0$ and hence \eq{Phisinh+} simplifies to \eq{bor} on the interval $]-\sinh(\lambda),\sinh(\lambda)[\,$. Using this, \eq{Phisinh+} can be rewritten as 
	\Eq{*}{
		d(u\sqrt{1+v^2}) 
		=\sqrt{1+v^2}d(u)+ud(\sqrt{1+v^2})
	}
	for all $u,v\in \,]-\sinh(\lambda),\sinh(\lambda)[\,$. With the substitution $w:=\sqrt{1+v^2}$, this equality yields that
	\Eq{*}{
		d(uw)=ud(w)+wd(u)
	}
	for all $u\in \,]-\sinh(\lambda),\sinh(\lambda)[\,$ and $w\in\,[1,\cosh(\lambda)[$. Then, in view of \lem{LRext}, $d$ is of Leibniz-type and hence is a standard derivation on $\R$.
	
	$(iii)$ In the third case suppose that $0<2\alpha<\beta$ (the other case can be treated similarly) and $d(\cosh(x))=\sinh(x)d(x)$ holds for all $x\in\,]\alpha,\beta[\,$. Then $\gamma=\alpha$ and $\delta=\frac\beta2$. The restriction of the cosine hyperbolic function to the interval $\,]\alpha,\beta[\,$ satisfies the functional equation 
	\Eq{*}{
		\cosh(x+y)=\cosh(x)\cosh(y)+\sqrt{\cosh^2(x)-1}\sqrt{\cosh^2(y)-1}
	}
	for all $x,y\in\,]\gamma,\delta[\,$. By \cor{add}, it follows that $d$ is also a derivation with respect to the function 
	\Eq{*}{
		(u,v)\mapsto uv+\sqrt{u^2-1}\sqrt{v^2-1}\qquad (u,v\in\,]\cosh(\gamma),\cosh(\delta)[\,),
	}
	i.e., the functional equation
	\Eq{Phicosh}{
		d(uv+\sqrt{u^2-1}\sqrt{v^2-1})=\Big(v+u\frac{\sqrt{v^2-1}}{\sqrt{u^2-1}}\Big)d(u)+\Big(u+v\frac{\sqrt{u^2-1}}{\sqrt{v^2-1}}\Big)d(v)
	}
	holds for all $u,v\in\,]\cosh(\gamma),\cosh(\delta)[\,$. Let $V$ be the set defined in \lem{Ratio}. Then, by this lemma, the intersection $V\cap \,]\cosh(\gamma),\cosh(\delta)[\,$ is nonempty, moreover it is dense in $]\cosh(\gamma),\cosh(\delta)[\,$. Then, for $u,v\in V\cap \,]\cosh(\gamma),\cosh(\delta)[\,$, we have that $u,v,$ and $\sqrt{v^2-1}\in\Q$. For such values of $u$ and $v$, the equality \eq{Phicosh} and the $\Q$-homogeneity of $d$ implies that
	\Eq{*}{
		(uv+\sqrt{u^2-1}&\sqrt{v^2-1})d(1)\\&
		=\Big(v+u\frac{\sqrt{v^2-1}}{\sqrt{u^2-1}}\Big)ud(1)+\Big(u+v\frac{\sqrt{u^2-1}}{\sqrt{v^2-1}}\Big)vd(1),
	}
	which is possible only if $d(1)=0$. Substituting $v:=u$ in \eq{Phicosh}, using the additivity and $\Q$-homogeneity of $d$ and that $d(1)=0$ we get that \eq{xr} with $r=2$ holds for all $u\in\,]\cosh(\gamma),\cosh(\delta)[\,$, which, using \lem{Nis}, implies that $d$ is a standard derivation.
	
	$(iv)$ In the fourth case assume that $\alpha<0<\beta$ and $d(\tanh(x))=\frac{1}{\cosh^2(x)}d(x)$ is valid for all $x\in\,]\alpha,\beta[\,$. Then $\gamma=\frac\alpha2$ and $\delta=\frac\beta2$. The restriction of the tangent hyperbolic function to the interval $\,]\alpha,\beta[\,$ satisfies the functional equation 
	\Eq{*}{
		\tanh(x+y)=\frac{\tanh(x)+\tanh(y)}{1+\tanh(x)\tanh(y)} \qquad(x,y\in\,]\gamma,\delta[\,).
	}
	By \cor{add}, it follows that $d$ is also a derivation with respect to the function 
	\Eq{*}{
		(u,v)\mapsto\frac{u+v}{1+uv}\qquad(u,v\in\,]\tanh(\gamma),\tanh(\delta)[\,),	
	}
	i.e., the functional equation
	\Eq{Phitanh}{
		d\Big(\frac{u+v}{1+uv}\Big)=\frac{1-v^2}{(1+uv)^2}d(u)+\frac{1-u^2}{(1+uv)^2}d(v)
	}
	holds for all $u,v\in\,]\tanh(\gamma),\tanh(\delta)[\,$. Now choose a subinterval $\,]\lambda,\mu[\,$ of $\,]\tanh(\gamma),\tanh(\delta)[\,$ such that $\lambda,\mu\in\Q$ and $0<\lambda\mu=:r$. It is easy to see that if $u\in \,]\lambda,\mu[\,$, then $\frac{r}{u}\in \,]\lambda,\mu[\,$ also holds. Substituting $u\in \,]\lambda,\mu[\,$ and $v:=\frac{r}{u}$, \eq{Phitanh} implies
	\Eq{*}{
		\frac1{1+r}d\Big(u+\frac{r}{u}\Big)=\frac{1-(\frac{r}{u})^2}{(1+r)^2}d(u)+\frac{1-u^2}{(1+r)^2}d\Big(\frac{r}{u}\Big).
	}
	A direct and simple computation yields that
	\Eq{*}{
		-\frac{1}{u^2}d(u)=d\Big(\frac{1}{u}\Big)
	}
	for $u\in \,]\lambda,\mu[\,$, which, using \lem{Nis}, implies that $d$ is a standard derivation.
	
	$(v)$ In the fifth case assume that $0<2\alpha<\beta$ (the other case can be handled similarly) and $d(\coth(x))=-\frac{1}{\sinh^2(x)}d(x)$ holds for all $x\in\,]\alpha,\beta[\,$. Then $\gamma=\frac\alpha2$ and $\delta=\frac\beta2$. The restriction of the cotangent hyperbolic function to the interval $\,]\alpha,\beta[\,$ satisfies the functional equation 
	\Eq{*}{
		\coth(x+y)=\frac{\coth(x)\coth(y)+1}{\coth(x)+\coth(y)} \qquad(x,y\in\,]\gamma,\delta[\,).
	}
	By \cor{add}, it follows that $d$ is also a derivation with respect to the function 
	\Eq{*}{
		(u,v)\mapsto\frac{uv+1}{u+v}\qquad(u,v\in\,]\coth(\delta),\coth(\gamma)[\,),	
	}
	i.e., the functional equation
	\Eq{Phicoth}{
		d\Big(\frac{uv+1}{u+v}\Big)=\frac{v^2-1}{(u+v)^2}d(u)+\frac{u^2-1}{(u+v)^2}d(v)
	}
	holds for all $u,v\in\,]\coth(\delta),\coth(\gamma)[\,$. Now choose a subinterval $\,]\lambda,\mu[\,$ of $\,]\coth(\delta),\coth(\gamma)[\,$ such that  $\lambda,\mu\in\Q$ and $\lambda+\mu=:r\neq0$. It is easy to see that if $u\in \,]\lambda,\mu[\,$, then $r-u\in \,]\lambda,\mu[\,$ also holds. Substituting $u\in \,]\lambda,\mu[\,$ and $v:=r-u$, \eq{Phicoth} implies
	\Eq{*}{
		\frac1r d(u(r-u)+1)=\frac{(r-u)^2-1}{r^2}d(u)+\frac{u^2-1}{r^2}d(r-u).
	}
	This equality, using the additivity and $\Q$-homogeneity, after some simplification, reduces to
	\Eq{*}{
		d(u^2)=2ud(u)+2d(1)-u^2d(1)\qquad(u\in\,]\lambda,\mu[\,).
	}
	If $u$ is rational, then this equality gives that $d(1)=0$. Hence the above equality shows that \eq{xr} is valid on $]\lambda,\mu[\,$ with $r=2$. In view of \lem{Nis}, we obtain that $d$ is a standard derivation. 
	
	$(vi)$ In the sixth case assume that $\alpha<0<\beta$ and $d(\sin(x))=\cos(x)d(x)$ holds for all $x\in\,]\alpha,\beta[\,$. Then $\gamma=\frac\alpha2$ and $\delta=\frac\beta2$. First we can choose $\lambda>0$ such that $\,]-\lambda,\lambda[\,\subseteq \,]\gamma,\delta[\,\cap\,]-\frac\pi2,\frac\pi2[\,$. Thus the cosine function is everywhere positive over $\,]-\lambda,\lambda[\,$. Then the sine function is strictly increasing on $\,]-\lambda,\lambda[\,$ and $\cos(x)=\sqrt{1-\sin^2(x)}$ holds for all $x\in\,]-\lambda,\lambda[\,$. Therefore the restriction of the sine function to the interval $\,]\alpha,\beta[\,$ satisfies the functional equation 
	\Eq{*}{
		\sin(x+y)=\sin(x)\sqrt{1-\sin^2(y)}+\sqrt{1-\sin^2(x)}\sin(y)
		\qquad (x,y\in\,]-\lambda,\lambda[\,).
	}
	By \cor{add}, it follows that $d$ is also a derivation with respect to the function 
	\Eq{*}{
		(u,v)\mapsto u\sqrt{1-v^2}+v\sqrt{1-u^2}\qquad (u,v\in\,]-\sin(\lambda),\sin(\lambda)[\,).
	}
	It means that the functional equation
	\Eq{Phisin}{
		d(u&\sqrt{1-v^2}+v\sqrt{1-u^2}) \\&
		=\Big(\sqrt{1-v^2}-v\frac{u}{\sqrt{1-u^2}}\Big)d(u)+\Big(\sqrt{1-u^2}-u\frac{v}{\sqrt{1-v^2}}\Big)d(v)
	}
	holds for all $u,v\in\,]-\sin(\lambda),\sin(\lambda)[\,$. Replacing $v$ by $-v$ and adding the equality so obtained to \eq{Phisin} side by side, we get that
	\Eq{Phisin+}{
		d(u\sqrt{1-v^2}) 
		=\sqrt{1-v^2}d(u)-u\frac{v}{\sqrt{1-v^2}}d(v)
	}
	holds for all $u,v\in\,]-\sin(\lambda),\sin(\lambda)[\,$. Let $W$ be the set defined in \lem{Ratio}. Then, by this lemma, the intersection $W\cap \,]-\sin(\lambda),\sin(\lambda)[\,$ is nonempty, moreover it is dense in $]-\sin(\lambda),\sin(\lambda)[\,$. Then, for $u,v\in W\cap \,]-\sin(\lambda),\sin(\lambda)[\,$, we have that $u,v,$ and $\sqrt{1-v^2}\in\Q$. For such values of $u$ and $v$, the equality \eq{Phisin+} and the $\Q$-homogeneity of $d$ implies that
	\Eq{*}{
		u\sqrt{1-v^2}d(1)
		=\sqrt{1-v^2}ud(1)-u\frac{v}{\sqrt{1-v^2}}vd(1),
	}
	which is possible only if $d(1)=0$. Let $0\neq u\in\,]-\sin(\lambda),\sin(\lambda)[\,\cap\,\Q$,
	and $v\in\,]-\sin(\lambda),\sin(\lambda)[\,$. By the $\Q$-homogeneity and additivity of $d$, it follows that $d(u)=ud(1)=0$ and hence \eq{Phisin+} simplifies to \eq{bor} on the interval $]-\sin(\lambda),\sin(\lambda)[\,$. Using this, \eq{Phisin+} can be rewritten as 
	\Eq{*}{
		d(u\sqrt{1-v^2}) 
		=\sqrt{1-v^2}d(u)+ud(\sqrt{1-v^2})
	}
	for all $u,v\in \,]-\sin(\lambda),\sin(\lambda)[\,$.
	With the substitution $w:=\sqrt{1-v^2}$, this equality yields that
	\Eq{*}{
		d(uw)=ud(w)+wd(u)
	}
	for all $u\in \,]-\sin(\lambda),\sin(\lambda)[\,$ and $w\in\,]\cos(\lambda),1]$. Then, in view of \lem{LRext}, $d$ is of Leibniz-type and hence is a standard derivation on $\R$.
	
	$(vii)$ In the seventh case suppose that $0<2\alpha<\pi<\beta$ (the other case can be treated similarly) and $d(\cos(x))=-\sin(x)d(x)$ holds for all $x\in\,]\alpha,\beta[\,$. Then $\gamma=\alpha$ and $\delta=\frac\beta2$ and we can choose $\lambda\in\,]0,\frac{\pi}{2}[\,$ such that $\,]\frac{\pi}{2}-\lambda,\frac{\pi}{2}+\lambda[\,\subseteq \,]\gamma,\delta[\,$. Thus the sine function is positive over $\,]\frac{\pi}{2}-\lambda,\frac{\pi}{2}+\lambda[\,$ and hence $\sin(x)=\sqrt{1-\cos^2(x)}$ holds for all $x\in\,]\frac{\pi}{2}-\lambda,\frac{\pi}{2}+\lambda[\,$. The restriction of the cosine function to the interval $\,]\alpha,\beta[\,$ satisfies the functional equation 
	\Eq{*}{
		\cos(x+y)=\cos(x)\cos(y)-\sqrt{1-\cos^2(x)}\sqrt{1-\cos^2(y)}
	}
	for all $x,y\in\,]\tfrac{\pi}{2}-\lambda,\tfrac{\pi}{2}+\lambda[\,$.
	By \cor{add}, it follows that $d$ is also a derivation with respect to the function 
	\Eq{*}{
		(u,v)\mapsto uv-\sqrt{1-u^2}\sqrt{1-v^2}\qquad (u,v\in\,]\cos(\tfrac{\pi}{2}+\lambda),\cos(\tfrac{\pi}{2}-\lambda)[\,),
	}
	i.e., the functional equation
	\Eq{Phicos}{
		d(uv-\sqrt{1-u^2}\sqrt{1-v^2})=\Big(v+u\frac{\sqrt{1-v^2}}{\sqrt{1-u^2}}\Big)d(u)+\Big(u+v\frac{\sqrt{1-u^2}}{\sqrt{1-v^2}}\Big)d(v)
	}
	holds for all $u,v\in\,]\cos(\tfrac{\pi}{2}+\lambda),\cos(\tfrac{\pi}{2}-\lambda)[\,$. Replacing $v$ by $-v$ and subtracting the equality so obtained from \eq{Phicos}, we get that
	\Eq{Phicos+}{
		d(uv)=vd(u)+ud(v)
	}
	holds for all $u,v\in\,]\cos(\tfrac{\pi}{2}+\lambda),\cos(\tfrac{\pi}{2}-\lambda)[\,$.
	Hence, in view of \lem{LRext}, we get that $d$ is a standard derivation on $\R$.
	
	$(viii)$ In the eighth case suppose that $-\pi<2\alpha<\beta$ and $\alpha<2\beta<\pi$ and $d(\tan(x))=\frac{1}{\cos^2(x)}d(x)$ is valid for all $x\in\,]\alpha,\beta[\,$. The restriction of the tangent function to the interval $\,]\alpha,\beta[\,$ satisfies the functional equation 
	\Eq{*}{
		\tan(x+y)=\frac{\tan(x)+\tan(y)}{1-\tan(x)\tan(y)}
		\qquad(x,y\in\,]\gamma,\delta[\,).
	}
	By \cor{add}, it follows that $d$ is also a derivation with respect to the function 
	\Eq{*}{
		(u,v)\mapsto\frac{u+v}{1-uv}\qquad(u,v\in\,]\tan(\gamma),\tan(\delta)[\,),	
	}
	i.e., the functional equation
	\Eq{Phitan}{
		d\Big(\frac{u+v}{1-uv}\Big)=\frac{1+v^2}{(1-uv)^2}d(u)+\frac{1+u^2}{(1-uv)^2}d(v)
	}
	holds for all $u,v\in\,]\tan(\gamma),\tan(\delta)[\,$. Now choose a subinterval $\,]\lambda,\mu[\,$ of $\,]\tan(\gamma),\tan(\delta)[\,$ such that $\lambda,\mu\in\Q$ and $0<\lambda\mu=:r\neq1$. It is easy to see that if $u\in \,]\lambda,\mu[\,$, then $\frac{r}{u}\in \,]\lambda,\mu[\,$ also holds. Substituting $u\in \,]\lambda,\mu[\,$ and $v:=\frac{r}{u}$, \eq{Phitan} implies
	\Eq{*}{
		\frac1{1-r}d\Big(u+\frac{r}{u}\Big)=\frac{1+(\frac{r}{u})^2}{(1-r)^2}d(u)+\frac{1+u^2}{(1-r)^2}d\Big(\frac{r}{u}\Big).
	}
	A direct and simple computation yields that
	\Eq{*}{
		-\frac{1}{u^2}d(u)=d\Big(\frac{1}{u}\Big)
	}
	for $u\in \,]\lambda,\mu[\,$, which, using \lem{Nis}, yields that $d$ is a standard derivation.
	
	$(ix)$ In the last case suppose that $0<2\alpha<\beta<\pi$ (the case $-\pi<\alpha<2\beta<0$ can be treated similarly) and $d(\cot(x))=-\frac{1}{\sin^2(x)}d(x)$ is valid for all $x\in\,]\alpha,\beta[\,$. The restriction of the cotangent function to the interval $\,]\alpha,\beta[\,$ satisfies the functional equation 
	\Eq{*}{
		\cot(x+y)=\frac{\cot(x)\cot(y)-1}{\cot(x)+\cot(y)} \qquad(x,y\in\,]\gamma,\delta[\,).
	}
	By \cor{add}, it follows that $d$ is also a derivation with respect to the function 
	\Eq{*}{
		(u,v)\mapsto\frac{uv-1}{u+v}\qquad(u,v\in\,]\cot(\delta),\cot(\gamma)[\,),	
	}
	i.e., the functional equation
	\Eq{Phicot}{
		d\Big(\frac{uv-1}{u+v}\Big)=\frac{1+v^2}{(u+v)^2}d(u)+\frac{1+u^2}{(u+v)^2}d(v)
	}
	holds for all $u,v\in\,]\cot(\delta),\cot(\gamma)[\,$. Now choose a subinterval $\,]\lambda,\mu[\,$ of $\,]\cot(\delta),\cot(\gamma)[\,$ such that  $\lambda,\mu\in\Q$ and $\lambda+\mu=:r\neq0$. It is easy to see that if $u\in \,]\lambda,\mu[\,$, then $r-u\in \,]\lambda,\mu[\,$ also holds. Substituting $u\in \,]\lambda,\mu[\,$ and $v:=r-u$, \eq{Phicot} implies
	\Eq{*}{
		\frac1r d(u(r-u)-1)=\frac{1+(r-u)^2}{r^2}d(u)+\frac{1+u^2}{r^2}d(r-u).
	}
	This equality, using the additivity and $\Q$-homogeneity, after some simplification, reduces to
	\Eq{*}{
		d(u^2)=2ud(u)-2d(1)-u^2d(1)\qquad(u\in\,]\lambda,\mu[\,).
	}
	If $u$ is rational, then this equality gives that $d(1)=0$. Hence the above equality shows that \eq{xr} is valid on $]\lambda,\mu[\,$ with $r=2$. In view of \lem{Nis} we get that $d$ is a standard derivation. 
\end{proof}

When considering the interval $\,]\alpha,\beta[\,$, one should observe that depending on additional assumptions described in the nine cases of the theorem, this interval can be arbitrary small in all the cases except the case $(vii)$, then it contains either $[\pi/2,\pi]$ or $[-\pi,-\pi/2]$.

\section{Extension of the result of Nishiyama and Horinouchi}

\Thm{PQ}{Let $d:\R\to\R$ be an additive function, $I$ be a nonempty open interval not containing zero and assume that $P,Q:I\to\R$ are of the form
\Eq{PQ0}{
  P(u)=\sum_{k\in\Z} p_ku^k,\qquad Q(u)=\sum_{k\in\Z} q_ku^k,
}
where $p_k,q_k\in\Q$ for all $k\in\Z$ and the set $\{k\in\Z\mid (p_k,q_k)\neq(0,0)\}$ is finite. Then 
\Eq{PQ}{
   Q'(u)d(P(u))=P'(u)d(Q(u))\qquad (u\in I)
}
holds if and only if  
\begin{enumerate}[(i)]                                                                                                                                      \item either $P$ and $Q$ are linearly dependent,
\item or $P$ and $Q$ are linearly independent, $P-p_0$ and $Q-q_0$ are linearly dependent and $d(1)=0$,
\item or $P-p_0$ and $Q-q_0$ are linearly independent and $d$ is a standard derivation.
\end{enumerate}}

\begin{proof} First we prove the necessity of (i)--(iii). 
Assume that \eq{PQ} holds and $P$, $Q$ are linearly independent. Then the Wronskian $P$ and $Q$ is not identically zero on $I$, i.e., there exists $u_0\in I$ such that $P'(u_0)Q(u_0)\neq P(u_0)Q'(u_0)$. 

Substituting $u\in I\cap\Q$ we have that $P(u)$ and $Q(u)$ are rational numbers. Therefore, using the $\Q$-homogeneity of $d$, \eq{PQ} implies that 
\Eq{PQu}{
  Q'(u)P(u)d(1)=P'(u)Q(u)ud(1)
}
for all $u\in I\cap\Q$. By the continuity of $P$ and $Q$ and the density of $I\cap\Q$ in $I$, it follows that \eq{PQu} holds for all $u\in I$, in particular, for $u=u_0$, which implies that $d(1)=0$. 

From now on, we assume that $P-p_0$ and $Q-q_0$ are linearly independent.
Let $k_0$ denote the smallest element of the finite set $\{k\in\Z\setminus\{0\}\mid (p_k,q_k)\neq(0,0)\}$. Then $u^{-k_0}(P(u)-p_0)$ and $u^{-k_0}(Q(u)-q_0)$ are linearly independent polynomials of the variable $u$. This is equivalent to the linear independence of the coefficients of $P-p_0$ and $Q-q_0$, that is, of $(p_k)_{k\in\Z\setminus\{0\}}$ and $(q_k)_{k\in\Z\setminus\{0\}}$. Therefore, the system of vectors $(p_k,q_k)_{k\in\Z\setminus\{0\}}$ spans $\R^2$. Thus there exists $\ell\in\Z\setminus\{0\}$ such that, for $0\neq k<\ell$, the vector $(p_k,q_k)$ is parallel to $(p_{k_0},q_{k_0})$ and $(p_\ell,q_\ell)$ is not  parallel to $(p_{k_0},q_{k_0})$. Then, 
\Eq{kl}{
  p_iq_j=p_jq_i \qquad (0\neq i<\ell,\,0\neq j<\ell) 
  \qquad\mbox{and}\qquad
  p_{k_0}q_\ell\neq p_\ell q_{k_0}.
}

Now let $v\in I$ be fixed and $r\in (I/v)\cap\Q)$. Then $u=rv\in I$, therefore, the equality \eq{PQ}, the additivity and $\Q$-homogeneity of $d$ imply
\Eq{*}{
  \bigg(\sum_{j\in\Z}q_jjr^{j-1}v^{j-1}\bigg)
  \bigg(\sum_{i\in\Z} p_ir^{i}d(v^i)\bigg)
  =\bigg(\sum_{i\in\Z}p_iir^{i-1}v^{i-1}\bigg)
  \bigg(\sum_{j\in\Z} q_jr^{j}d(v^j)\bigg).
}
Using that $d(1)=0$, this equality is equivalent to
\Eq{ij}{
  \sum_{i\in\Z\setminus\{0\}}\sum_{j\in\Z\setminus\{0\}} 
  p_iq_jr^{i+j-1}(jv^{j-1}d(v^i)-iv^{i-1}d(v^j))=0.
}
This implies
\Eq{*}{
  \sum_{i<j,\,ij\neq0} 
  (p_iq_j-p_jq_i)r^{i+j-1}(jv^{j-1}d(v^i)-iv^{i-1}d(v^j))=0.
}
According to the choice of $k_0$ and $\ell$, we have \eq{kl}, therefore,
\Eq{*}{
  \sum_{k_0\leq i,\,\max(i+1,\ell)\leq j,\,ij\neq0} 
  (p_iq_j-p_jq_i)r^{i-k_0+j-\ell}(jv^{j-1}d(v^i)-iv^{i-1}d(v^j))=0.
}
The left hand side of this equality is a polynomial of $r$, hence its value at $r=0$ is equal to zero, which gives 
\Eq{*}{
  (p_{k_0}q_\ell-p_\ell q_{k_0})(\ell v^{\ell-1}d(v^{k_0})-k_0v^{k_0-1}d(v^\ell))=0.
}
By the last relation in \eq{kl}, this yields
\Eq{*}{
  \ell d(v^{k_0})=k_0v^{k_0-\ell}d(v^\ell) \qquad(v\in I).
}
With the substitution $u:=v^\ell$, and with the notation $r:=\frac{k_0}{\ell}$, we get
\Eq{*}{
  d(u^r)=ru^{r-1}d(u) \qquad(u\in J:=\{x^\ell\mid x\in I\}).
}
Observe that $r\in\Q\setminus\{0,1\}$, therefore, by \lem{Nis}, it follows that $d$ is a standard derivation.

If condition (i) holds and $P$ is not identically zero, then $Q$ is a rational multiple of $P$, hence \eq{PQ} is trivially valid by the $\Q$-homogeneity of $d$.

If condition (ii) holds, then, denoting $P_0:=P-p_0$ and $Q_0:=Q-q_0$ and using (i) for $P_0$ and $Q_0$, we have
\Eq{*}{
  Q'(u)&d(P(u))=Q_0'(u)d(P_0(u)+p_0)
  =Q_0'(u)(d(P_0(u))+p_0d(1))=Q_0'(u)d(P_0(u))\\
  &=P_0'(u)d(Q_0(u))=P_0'(u)(d(Q_0(u))+q_0d(1))
  =P_0'(u)d(Q_0(u)+q_0)=P'(u)d(Q(u)).
}
Finally, if condition (iii) is valid, i.e., $d$ is a standard derivation, then, for all $u\in I$ and $i,j\in\Z$,
\Eq{*}{
   ju^{j-1}d(u^i)=ju^{j-1}iu^{i-1}d(u)=iu^{i-1}d(u^j).
}
Multiplying this equality by $p_iq_j$, then summing up the equalities so obtained side by side for $(i,j)\in\Z^2$, we obtain that 
\Eq{*}{
  \bigg(\sum_{j\in\Z}q_jju^{j-1}\bigg)
  \bigg(\sum_{i\in\Z} p_id(u^i)\bigg)
  =\bigg(\sum_{i\in\Z}p_iiu^{i-1}\bigg)
  \bigg(\sum_{j\in\Z} q_jd(u^j)\bigg).
}
This equality, by the additivity and $\Q$-homogeneity of $d$ is equivalent to \eq{PQ}.
\end{proof}

\Cor{PQ}{Let $d:\R\to\R$ be an additive function, $I$ be a nonempty open interval not containing zero and assume that $P,Q:I\to\R$ are of the form \eq{PQ0}, where $p_k,q_k\in\Q$ for all $k\in\Z$ and the set $\{k\in\Z\mid (p_k,q_k)\neq(0,0)\}$ is finite. Assume that $Q'$ is non-vanishing on $I$, furthermore $P-p_0$ and $Q-q_0$ are linearly independent and $d$ derivates $P\circ Q^{-1}$ on $J:=Q(I)$. Then $d$ is a standard derivation.}

\begin{proof}
Using that $d$ derivates $P\circ Q^{-1}$, we have  
\Eq{*}{
	d(P(Q^{-1}(v)))=P'(Q^{-1}(v))\frac{1}{Q'(Q^{-1}(v))}d(v)\qquad(v\in J).
} 
Substituting $u:=Q^{-1}(v)\in I$ into the above equation, we can see that \eq{PQ} holds. Taking into consideration that $P-p_0$ and $Q-q_0$ are linearly independent by assumption, \thm{PQ} implies that $d$ is a standard derivation.
\end{proof}

\Cor{pol}{Let $d:\R\to\R$ be an additive function, let $P$ be a nonzero real polynomial with rational coefficients and $I$ be a nonempty open subinterval of $\R$. Then $d$ derivates $P$ over $I$ if and only if 
\begin{enumerate}[(i)]                                                                                                                                      \item either $\deg(P)=0$ and $d(1)=0$,
\item or $\deg(P)=1$ and $P(0)d(1)=0$,
\item or $\deg(P)\geq2$ and $d$ is a standard derivation.
\end{enumerate}
}

\begin{proof}
 This statement is an immediate consequence of \thm{PQ} by choosing $Q(u)=u$.
\end{proof}

\section{Open Questions}

Motivated by the results of the previous section, we can formulate two open problems. Let $P,Q:I\to\R$ by of the form \eq{PQ0}, where $p_k,q_k\in\Q$ for all $k\in\Z$ and the set $\{k\in\Z\mid (p_k,q_k)\neq(0,0)\}$ is finite and let $d:\R\to\R$ be an additive function.

\textbf{Problem 1.} Assume that $Q$ is non-vanishing on $I$ and $d$ derivates $P/Q$ on $I$, that is,
\Eq{*}{
  d\bigg(\frac{P(u)}{Q(u)}\bigg)
  =\frac{P'(u)Q(u)-Q'(u)P(u)}{Q^2(u)}d(u) \qquad(u\in I).
}
Under what conditions on $P$ and $Q$ does this equality imply that $d$ is a standard derivation? 

\textbf{Problem 2.} Assume that $P'$ is non-vanishing on $I$ and $Q(I)\subseteq P(I)$ and $d$ derivates $P^{-1}\circ Q$ on $I$, that is,
\Eq{*}{
  d\big(P^{-1}(Q(u)\big)=\frac{Q'(u)}{P'(P^{-1}(Q(u)))} d(u) \qquad(u\in I).
}
This, provided that $Q$ is strictly monotone on $I$, is equivalent to the
condition
\Eq{*}{
  P'(P^{-1}(v))d\big(P^{-1}(v)\big)=Q'(Q^{-1}(v)) d(Q^{-1}(v)) \qquad(v\in Q(I)).
}
Under what conditions on $P$ and $Q$ does this equality imply that $d$ is a standard derivation? 

The result of Boros and Erdei \cite{BorErd05} would be a particular case of such a generalization.


\end{document}